\documentclass[pdflatex,sn-mathphys-num]{sn-jnl}

\usepackage[T1]{fontenc}
\usepackage{lmodern}
\usepackage{amsmath,amssymb,amsthm,mathtools}
\usepackage{graphicx,booktabs,pgfplots}
\usepackage{xcolor}
\hypersetup{hidelinks}

\pgfplotsset{compat=1.18}
\newcommand{\R}{\mathbb R}
\newcommand{\I}{\mathcal I}
\newcommand{\dd}{\,\mathrm d}
\newcommand{\diver}{\operatorname{div}}
\newcommand{\tr}{\operatorname{tr}}
\newcommand{\Hess}{\operatorname{Hess}}

\newtheorem{theorem}{Theorem}
\newtheorem{proposition}{Proposition}
\newtheorem{lemma}{Lemma}
\theoremstyle{definition}
\newtheorem{definition}{Definition}
\newtheorem{example}{Example}

\begin{document}

\title[Optimal volume contraction and regularity]{Optimal Weighted Volume Dissipativity and Sobolev Regularity at Limiting Rate}


\author*[1]{\fnm{I.B.} \sur{Furtat}}\email{cainenash@mail.ru}

\affil[1]{\orgname{IPME RAS}, \orgaddress{\city{St.~Petersburg}, \country{Russia}}}


\abstract{Sharp pointwise divergence bounds for a nonlinear flow are studied under the requirement that the logarithm of the adapted positive density also be an exponentially decaying Lyapunov function. Invariant-measure duality supplies the unrestricted benchmark. Under two-sided quadratic Lyapunov bounds, scaling a prescribed Lyapunov function realizes every strictly suboptimal pointwise rate without reducing its decay rate, and the closure of the joint rate region is a rectangle. For a fixed reference volume form, the limiting pointwise rate has a smooth compatible certificate exactly when the differential of the reference divergence vanishes at the equilibrium. Otherwise, the infimal second-derivative supremum norm grows inversely with the gap to the limiting rate. Matching Sobolev cost bounds are proved in every finite dimension when the integrability exponent exceeds the dimension and the prescribed Lyapunov rate does not exceed half the baseline rate. The upper construction uses a regularized square root of the baseline Lyapunov function and keeps the potential and its gradient uniformly bounded. Matching contraction deficits under a Sobolev norm budget follow. A scalar borderline family has bounded second-derivative mass and a Lipschitz limit with atomic second derivative. The limiting exponent also holds in the reference measure with an explicit uniform volume prefactor, even when no smooth compatible instantaneous certificate exists. Scaling preserves the full baseline Lyapunov rate.}



\maketitle

\section{Introduction}\label{sec:intro}
Phase-space volume describes the collective evolution of trajectories. Liouville's formula relates its change to the divergence of the vector field~\cite{Liouville}, and transport identities describe integrals over moving regions~\cite{Reynolds}. The geometric language of flows, volume forms and vector fields is developed in~\cite{Arnold,Krasnoselski63}. Changing the volume density can reveal contraction absent in the original coordinates, raising the question of the largest uniform contraction rate obtainable in this way.

Divergence and index methods have a substantial history. Early studies related divergence to differential equations~\cite{Zaremba54}, connected Liouville's theorem with the stability problem~\cite{Fronteau65}, and investigated index, divergence and stability in autonomous systems~\cite{Brauchli68}. Subsequent developments provided source-based stability methods~\cite{Jukov78,Jukov79}, instability conditions~\cite{Jukov90} and asymptotic stability criteria for second-order systems~\cite{Jukov99}. These approaches require hypotheses beyond the sign of divergence: volume contraction alone does not ensure stability of an individual trajectory. Volume dissipativity below refers to transported phase-space regions, rather than the storage--supply formulation~\cite{Willems72,Willems72II}.

Modern density methods make the connection with trajectory convergence more explicit. Density certificates for almost-everywhere convergence and a converse theorem were established in~\cite{Rantzer01,RantzerConverse2002}. Convex density inequalities support nonlinear control synthesis~\cite{PrajnaParriloRantzer2004}, while Lyapunov measures provide an operator-theoretic treatment~\cite{VaidyaMehta2008}. Combined Lyapunov and density conditions were developed for rotational motion~\cite{VasconcelosRantzer2009,VasconcelosRantzer2011}. Lyapunov-density criteria and converse results for time-varying and periodic systems~\cite{Masubuchi2020}, and extended conditions for robust stabilization~\cite{Masubuchi2023}, complement classical Lyapunov theory~\cite{Khalil2002}. Our additional constraint is specific: the density is smooth and positive at an equilibrium, and its logarithm itself must be an exponentially decaying Lyapunov function.

Further developments concern stability and control. Flow and divergence criteria for stability and control were obtained in~\cite{Furtat2020ARC}, nonautonomous stability and feedback design were considered in~\cite{FurtatGushchin2020JFI}. Autonomous stability and instability conditions~\cite{FurtatGushchin2021Access} were followed by divergence criteria for exponential stability~\cite{FurtatGushchin2022Access}. Density systems were studied through quadratic functions whose derivatives depend on a density~\cite{Furtat2023Density}, and adaptive control of perturbed density systems was developed in~\cite{Furtat2026TAC}. The present problem fixes the autonomous vector field and optimizes the reference-volume adaptation, subject to the logarithm of its positive density being a Lyapunov function. Its specific questions are the joint optimal values, endpoint attainment and derivative costs.

Weighted volume and conformal geometry are considered in~\cite{Furtat26_Weighted}. A finite-dimensional log-affine family for regional phase-volume contraction is studied in~\cite{Furtat26_Density}. A family-restricted optimum need not identify the unrestricted smooth value. For the latter, established invariant-measure duality for extremal time averages~\cite{TobascoGoluskinDoering2018} supplies a benchmark when applied to divergence. The coboundary viewpoint is also standard in ergodic optimization~\cite{Jenkinson2019}. Neither framework alone imposes quadratic bounds, a minimum at the equilibrium and Lyapunov decay on the auxiliary potential. Differential contraction methods~\cite{Lohmiller98,Simpson14,Forni14} instead control tangent directions. The scalar volume objective measures their total expansion or contraction.

In~\cite{GarciaLatorreMaza2026}, inverse Jacobi multipliers are
used to locate periodic orbits and study their orbital stability.
Where a multiplier is nonzero, the reciprocal of its absolute
value gives a positive density with respect to which the extended
vector field has zero divergence. The time-dependent linear
changes of variables studied in~\cite{GaetaWalcher2025} transform
an autonomous equation into a nonautonomous one. In our problem,
the autonomous vector field is fixed. We seek a positive density
whose logarithm is a Lyapunov function for the equilibrium and
determine how close the weighted divergence can come to its
optimal uniform bound under this requirement.
Bohl, Lyapunov and Sacker--Sell spectra for nonautonomous linear equations are compared in~\cite{Doan2017}. Series conditions characterizing asymptotic stability and instability, together with perturbation results, are developed in~\cite{Dragicevic2025}. Nonuniform hyperbolicity is characterized by admissibility between weighted function spaces in~\cite{Chemnitz2025}. Entropy spectra of Lyapunov exponents are related to pressure for typical cocycles in~\cite{Mohammadpour2025}. Here the variational objective concerns the determinant cocycle, and the function-space question concerns regularity of scalar density potentials approaching its extremal divergence bound.

We ask whether one potential can certify both trajectory decay and nearly optimal volume contraction on a prescribed compact positively invariant region. Does compatibility lower the unrestricted value? Is the boundary attained, and what derivative cost is unavoidable when it is not? The central results assume an interior equilibrium and a smooth Lyapunov function with two-sided quadratic bounds and exponential decay on the entire region. This hypothesis forces the sole invariant probability to be the equilibrium mass. Although the optimal value is determined by the equilibrium trace, its realization by compatible densities depends on nonlinear regional constraints and may require unbounded second derivatives. The main contribution is an exact characterization of joint rate feasibility, smooth boundary attainment, and the sharp regularity cost of approaching the optimal divergence bound.

The contributions are as follows.
\begin{itemize}
\item Section~\ref{sec:joint} proves that scaling a prescribed baseline Lyapunov function realizes every strictly suboptimal volume rate without reducing its decay rate. It determines the rectangular closure of the joint rate region and its sharp local spectral bounds. This quantifies simultaneous optimization with a single logarithmic potential beyond density-existence and combined Lyapunov--density criteria~\cite{RantzerConverse2002,VasconcelosRantzer2011,Masubuchi2020}.
\item Section~\ref{sec:endpoint} gives an exact smooth endpoint criterion relative to the fixed reference form: the differential of its divergence must vanish at the equilibrium. Otherwise the infimal Hessian supremum norm grows inversely with the gap to the optimal rate. An additional integral volume estimate realizes the limiting exponent in the reference measure with an explicit uniform prefactor, including when smooth compatible pointwise attainment fails.
\item Section~\ref{sec:sobolev} determines the exact asymptotic Sobolev cost of approaching the optimal contraction rate when the integrability exponent exceeds the space dimension and the prescribed Lyapunov rate does not exceed half the baseline rate. A lower bound valid for every admissible potential is matched by an explicit regularized square-root construction with uniformly bounded density and first derivatives. The resulting sharp contraction deficits under a norm budget quantify the regularity required for near-optimal contraction. A scalar borderline example exhibits concentration of second derivatives into an atomic measure.
\end{itemize}

Section~\ref{sec:dual} supplies the unrestricted benchmark and an explicit example of its nonattainment. Section~\ref{sec:computation} treats constrained searches, and Section~\ref{sec:examples} gives examples. Appendix~\ref{app:transport} contains transport and conditional stability consequences. Appendix~\ref{app:perturbation} retains an alternative compatibility construction and a static sign test. The reference form remains fixed in the compatibility problem: changing it preserves the unrestricted value but changes the admissible density logarithms.

\section{Volume forms and the unrestricted value}\label{sec:dual}
Let $f$ be a smooth vector field on a neighborhood $U$ of a nonempty compact positively invariant set $K$ in an oriented $n$-dimensional manifold, $n\ge1$. Its forward flow $\varphi_t$ is defined on $K$ for all $t\ge0$. A compact manifold with inward-pointing field at its boundary is allowed. Fix a smooth positive volume form $m$, for instance $dV_g$ for a background metric $g$. Write
\[
 b=\diver_m f,\qquad f(u)=du[f],\qquad D_u=b+f(u).
\]
Here $du$ is the differential of $u$, so $f(u)$ is its derivative along the vector field, and $\mathcal L_f$ is the Lie derivative. We use $C(K)$ for continuous real functions with the uniform norm, $C^k(U)$ for $k$-times continuously differentiable functions, and $C_c^\infty(U)$ for smooth functions with compact support. Smoothness at a boundary is understood up to that boundary. In fixed coordinates, $D^2u$ is the Hessian matrix, $Df(p)$ is the linearization at an equilibrium, $\top$ denotes transpose and $I$ the identity matrix. Matrix inequalities refer to positive semidefiniteness, $H\succ0$ means positive definiteness. A matrix is Hurwitz when its spectrum lies in the open left half-plane. We write $(z)_+=\max\{z,0\}$, $E\Subset U$ for compact containment, and $\overline E$ and $\operatorname{int}E$ for closure and interior. The symbols $O$ and $o$ refer to bounded and vanishing ratios in the specified limit.

Divergence is defined by $\mathcal L_fm=bm$, so
\begin{equation}\label{eq:change}
 \diver_{e^um}f=b+f(u).
\end{equation}
For $m=dV_g$, the conformal metric $e^{2u/n}g$ has volume form $e^um$. If $dV_h=\rho dV_g$, the metric $\rho^{2/n}g$ has exactly the same volume and divergence as $h$. Thus unrestricted searches over densities, conformal metrics and all metrics have the same value. This does not concern metric classes with curvature, boundary or normalization constraints.

Positive invariance means $\varphi_t(K)\subset K$ for all $t\ge0$. Let $\I_f(K)$ denote invariant Borel probability measures supported in $K$: $\int v\circ\varphi_t\,d\mu=\int v\,d\mu$ for every $v\in C(K)$. Weak-* convergence means convergence against every such $v$, $\delta_x$ is the Dirac mass at $x$. Define
\begin{equation}\label{eq:values}
 \beta_K(f)=\max_{\mu\in\I_f(K)}\int_Kb\dd\mu,
 \qquad \lambda_*=
 \sup_{u\in C^\infty(U)}\inf_K[-b-f(u)].
\end{equation}
The latter value may be nonpositive. A pointwise bound $D_u\le-\lambda$ implies, for measurable $E\subset K$,
\begin{equation}\label{eq:volume}
 \int_{\varphi_t E}e^u m\le e^{-\lambda t}\int_Ee^um,
 \qquad t\ge0,
\end{equation}
by transport (Appendix~\ref{app:transport}).

\begin{proposition}\label{th:dual}
The values satisfy
\begin{equation}\label{eq:duality}
 \beta_K(f)=\inf_{u\in C^\infty(U)}\max_K(b+f(u)),
 \qquad \lambda_*=-\beta_K(f),
\end{equation}
and are independent of the reference form and $U$. A smooth positive density, or a smooth metric, with strictly negative divergence on $K$ exists exactly when $\beta_K(f)<0$. The finite-time potentials
\begin{equation}\label{eq:average}
 u_T(x)=\frac1T\int_0^T(T-t)b(\varphi_t x)\dd t
\end{equation}
satisfy
\begin{equation}\label{eq:averageidentity}
 b+f(u_T)=q_T:=\frac1T\int_0^Tb(\varphi_t x)\dd t,
 \qquad \max_Kq_T\rightarrow\beta_K(f).
\end{equation}
They extend smoothly to $U$ without changing their jets on $K$.

The spectral identity below requires only $C^1$ regularity of $f$. If $K$ is fully invariant, then
\begin{equation}\label{eq:exponents}
 \beta_K(f)=\max_{\mu\in\I_f^{\rm erg}(K)}\sum_{j=1}^n\chi_j(\mu),
\end{equation}
where $\I_f^{\rm erg}(K)$ consists of ergodic invariant probabilities and the tangent-flow Lyapunov exponents $\chi_j(\mu)$ are counted with multiplicity. If periodic orbit measures are weak-* dense in $\I_f(K)$, then
\begin{equation}\label{eq:periodic}
 \beta_K(f)=\sup_{\gamma\subset K}\frac1{T_\gamma}
 \log|\det D\varphi_{T_\gamma}(x_\gamma)|.
\end{equation}
Here the supremum is over periodic orbits $\gamma$, with period $T_\gamma>0$ and $x_\gamma\in\gamma$, and $D\varphi_t$ is the differential of the flow. If $f_j\to f$ in $C^1$ near a fixed common positively invariant compact set $K$, then $\limsup_j\beta_K(f_j)\le\beta_K(f)$. If $K$ is an entire compact manifold without boundary, then $\beta_K(f)\ge0$.
\end{proposition}
\begin{proof}
Empirical measures on compact $K$ have weak-* accumulation points, and limits for intervals tending to infinity are invariant.
Hence $\I_f(K)$ is nonempty and compact. For invariant $\mu$, differentiation of $\int u\circ\varphi_t\dd\mu$ gives $\int f(u)\dd\mu=0$, so $\int b\dd\mu\le\max_K(b+f(u))$.
Integration by parts gives~\eqref{eq:averageidentity}. Invariance gives $\max_Kq_T\ge\beta_K(f)$. For the reverse bound choose maximizing points $x_j$ for times $T_j\to\infty$. Every limit of
$T_j^{-1}\int_0^{T_j}\delta_{\varphi_t x_j}\dd t$
is invariant: a shift by fixed $s$ changes its integral against $v\in C(K)$ by at most $2s\|v\|_\infty/T_j$. Thus every limiting maximal average is at most $\beta_K(f)$. This proves the convergence and duality. For fixed $T$, compactness of $K$ and smooth dependence on initial conditions give an open neighborhood $W_T\subset U$ on which the flow exists up to time $T$ and $u_T$ is smooth. Choose a smooth cutoff supported in $W_T$ and equal to one on a neighborhood of $K$. Multiplying $u_T$ by this cutoff and extending by zero gives a smooth function on $U$, unchanged near $K$. Changing $m$ to $e^qm$ adds $f(q)$ to $b$, which has zero invariant average. This proves independence and the negative-divergence equivalence.

Writing $J_m\varphi_t$ for the positive flow Jacobian relative to $m$, Liouville's formula gives
\begin{equation}\label{eq:jacobian}
 \log J_m\varphi_t(x)=\int_0^tb(\varphi_sx)\dd s.
\end{equation}
Birkhoff's theorem and the multiplicative ergodic theorem~\cite{Oseledets1968} identify its almost-sure time average under ergodic $\mu$ with both $\int b\dd\mu$ and $\sum_j\chi_j(\mu)$. A maximizing measure has a maximizing ergodic component, proving~\eqref{eq:exponents}. Periodic-measure density and~\eqref{eq:jacobian} give~\eqref{eq:periodic}. Density is an additional hypothesis, available for suitable hyperbolic flows~\cite{Sigmund1972}, not an automatic property.

For robustness take maximizing $\mu_j$ along a subsequence realizing the upper limit. A weak-* subsequential limit $\mu$ is invariant under $f$, by uniform finite-time convergence of the flows. Uniform convergence of the divergences gives
$\lim_j\int\diver_m f_j\dd\mu_j=\int b\dd\mu\le\beta_K(f)$.
Thus strict negativity is robust for a fixed $K$. Continuation of varying invariant sets needs an additional hypothesis. Finally, on a closed manifold Stokes gives $\int\diver_{\hat m}f\,\hat m=0$ for every positive smooth form $\hat m$. Uniformly negative divergence is impossible, implying $\beta_K(f)\ge0$.
\end{proof}
The unrestricted value need not be attained, even when it is positive, as the following example proves. Building on auxiliary-function duality~\cite{TobascoGoluskinDoering2018}, Proposition~\ref{th:dual} characterizes the unrestricted divergence optimum and supplies explicit nearly optimal potentials. The sum in~\eqref{eq:exponents} contains every exponent. No entropy identity or absolute continuity hypothesis is involved. In particular, it differs from the entropy formula involving positive exponents under the additional hypotheses of smooth ergodic theory~\cite{Pesin1977}. The Stokes obstruction permits searches on proper invariant subsets, spaces with boundary and noncompact spaces, but these settings alone do not ensure a negative value.

\begin{example}\label{ex:unrestricted}
On $K=[-1,1]^2$ with Euclidean volume, let
\[
 f(x,y)=(-x^3,(-1+ax)y),\qquad 0<a<1.
\]
The field points inward on each face. Invariance and $f(x^2)=-2x^4$ force every invariant probability to be supported on $x=0$. Then $f(y^2)=-2y^2$ forces it to be $\delta_{(0,0)}$. Since $b=-1+ax-3x^2$, Proposition~\ref{th:dual} gives $\lambda_*=1$. For any smooth potential $u$,
\[
 b(x,0)+f(u)(x,0)+1=ax-3x^2-x^3\partial_xu(x,0)>0
\]
for all sufficiently small $x>0$, because $\partial_xu$ is bounded near zero. Thus no smooth density satisfies $D_u\le-1$ on $K$, although every smaller positive margin is feasible. Indeed, the first coordinate satisfies $x(t)=x(0)/\sqrt{1+2t x(0)^2}$, so the averages in~\eqref{eq:averageidentity} obey $q_T\le-1+a\sqrt{2/T}$ uniformly on $K$. This is nonattainment in the unrestricted density class. The linearization has a zero eigenvalue, so the quadratic exponential hypothesis below does not apply.
\end{example}

\section{Compatibility without loss of the optimal value}\label{sec:joint}
Let $p\in\operatorname{int}K$ be an equilibrium, $d=d(x,p)$ denotes distance in a fixed smooth background metric. Assume some $V\in C^\infty(U)$ satisfies
\begin{equation}\label{eq:base}
 V(p)=0,\quad a_0d^2\le V\le A_0d^2,\quad
 f(V)\le-c_0V\quad(x\in K),\qquad a_0,A_0,c_0>0.
\end{equation}
The quadratic constants are not fixed during optimization.
\begin{definition}\label{def:compatible}
A smooth $u$ is compatible with rates $(\lambda,c)\in(0,\infty)^2$ if, for some $a,A>0$,
\begin{equation}\label{eq:compatible}
 u(p)=0,\quad ad^2\le u\le Ad^2,\quad
 f(u)\le-cu,\quad b+f(u)\le-\lambda\quad(x\in K).
\end{equation}
Let $\mathcal J_K$ be the set of such pairs. Let $c_L(K)$ be the supremum of positive decay rates under the first three conditions, with value zero for an empty class.
\end{definition}
Integration gives the simultaneous estimates
\begin{equation}\label{eq:distance}
 \begin{split}
 u(\varphi_tx)&\le e^{-ct}u(x),\\
 d(\varphi_tx,p)&\le\sqrt{A/a}\,e^{-ct/2}d(x,p),
 \end{split}
\end{equation}
Together with~\eqref{eq:volume}, these give both estimates. Thus $c$ is decay of the potential and $c/2$ is the guaranteed distance rate. The density $e^u$ remains positive at $p$, with value one. Related conditional stability criteria for vanishing weights, invariant submanifolds and strongly convex gradient potentials are given in Appendix~\ref{app:transport}. They require the additional hypotheses stated there.

For local statements fix a coordinate neighborhood $U_0$ of $p$. Let $\mathfrak K_p$ consist of the compact positively invariant neighborhoods of $p$ contained in $U_0$, and set
\[
 \mathcal J_{\rm loc}=\bigcup_{K\in\mathfrak K_p}\mathcal J_K.
\]
Define $c_L^{\rm loc}(p)$ as the supremum of positive $c$ for which a $C^2$ function on some neighborhood of $p$ satisfies $u(p)=0$, two-sided quadratic bounds and $f(u)\le-cu$. Take zero if this class is empty. Neither local definition presupposes~\eqref{eq:base}.

\begin{theorem}\label{th:joint}
Under~\eqref{eq:base}, write $A_f=Df(p)$ and $\alpha(A_f)=\max\{\operatorname{Re}z:z\in\sigma(A_f)\}$, where $\sigma$ denotes spectrum and $\operatorname{Re}$ real part. Then $A_f$ is Hurwitz and
\begin{equation}\label{eq:unique}
 \I_f(K)=\{\delta_p\},\qquad \lambda_*=-b(p)=-\tr A_f>0.
\end{equation}
For every $0<\lambda<\lambda_*$,
\begin{equation}\label{eq:scaling}
 S_\lambda=\sup_{x\in K\setminus\{p\}}
 \frac{(b(x)+\lambda)_+}{-f(V)(x)}<\infty,
\end{equation}
and every $s>0$, $s\ge S_\lambda$, gives a compatible $u=sV$ with rates $(\lambda,c_0)$. Consequently,
\begin{equation}\label{eq:rectangle}
 \begin{split}
 (0,\lambda_*)\times(0,c_L(K))&\subset\mathcal J_K
 \subset(0,\lambda_*]\times(0,c_L(K)],\\
 \overline{\mathcal J_K}&=[0,\lambda_*]\times[0,c_L(K)].
 \end{split}
\end{equation}
Independently of the regional assumption~\eqref{eq:base}, local compatibility on some compact positively invariant neighborhood, local exponential stability, and the Hurwitz property of $A_f$ are equivalent. Under these equivalent local conditions,
\begin{equation}\label{eq:localrectangle}
 \overline{\mathcal J_{\rm loc}}=[0,-\tr A_f]\times[0,-2\alpha(A_f)].
\end{equation}
The volume value equals $-\tr A_f$ on every sufficiently small Lyapunov sublevel neighborhood. For a linear Hurwitz system, $c_L(K)=-2\alpha(A_f)$ on every compact positively invariant ellipsoid containing the origin in its interior.
\end{theorem}
\begin{proof}
Integration of $f(V)\le-c_0V$ against invariant $\mu$ gives $\int V\dd\mu=0$, hence $\mu=\delta_p$. Proposition~\ref{th:dual} gives $\lambda_*=-b(p)$. In coordinates at $p$, the quadratic bounds imply $H=D^2V(p)\succ0$, and the second-order decay inequality gives
\begin{equation}\label{eq:matrix}
 A_f^\top H+HA_f\preceq-c_0H.
\end{equation}
Thus $A_f$ is Hurwitz, at an equilibrium divergence relative to any smooth positive form equals $\tr A_f$. This proves~\eqref{eq:unique}.

For $\lambda<\lambda_*$, the numerator in~\eqref{eq:scaling} vanishes near $p$. On the compact remainder, $-f(V)\ge c_0V$ has a positive minimum. Hence $S_\lambda$ is finite. The choice of $s$ gives $b+sf(V)\le-\lambda$ and $f(sV)\le-c_0sV$. Scaling preserves $A_0/a_0$. For any $c<c_L(K)$ choose a $V$ with rate $c_0>c$ and apply this construction. Evaluation at $p$ gives $\lambda\le-\tr A_f$, the definition of $c_L$ gives the other upper bound. Taking closures proves~\eqref{eq:rectangle}.

The sharp local quadratic Lyapunov rate is $c_L^{\rm loc}(p)=-2\alpha(A_f)$. Indeed, a $C^2$ potential with rate $c$ yields~\eqref{eq:matrix} with $c$ replacing $c_0$. Testing on a complex eigenvector gives $c\le-2\alpha(A_f)$. Conversely, for $0<c<-2\alpha(A_f)$ choose $P\succ0$ with
$A_f^\top P+PA_f=-cP-I$.
The local function $V(x)=x^\top Px$ satisfies $f(V)=-cV-|x|^2+o(|x|^2)\le-cV$ after shrinking the coordinate neighborhood. For sufficiently small $\rho>0$, the sublevel $K_\rho=\{V\le\rho\}$ lies inside that neighborhood and is compact. Since $f(V)<0$ on its boundary, it is positively invariant. Quadratic bounds and $f(V)\le-cV$ give local exponential stability. Conversely, differentiating a local exponential flow estimate at $p$ yields $\|e^{tA_f}\|\le M e^{-\omega t}$ for some $M,\omega>0$, so $A_f$ is Hurwitz. Compatibility implies stability by~\eqref{eq:distance}.

Apply scaling on these sublevels, choosing a local Lyapunov rate strictly above the requested $c$, to obtain the interior of~\eqref{eq:localrectangle}. The equilibrium and eigenvector tests give its upper bounds. In the linear case the quadratic construction applies on any such ellipsoid, proving the regional formula as well.
\end{proof}
Appendix~\ref{app:perturbation} gives an alternative proof of compatibility at strictly suboptimal rates: Lemma~\ref{lem:jet} corrects the differential of a selected nearly optimal auxiliary potential at $p$, after which adding $sV$ gives Lyapunov decay rates tending to $c_0$. This recovers the joint rectangle while retaining that auxiliary potential in the construction.

There is no trade-off between the two supremal values in the freely scalable smooth class. Fixed bases, amplitudes, quadratic constants or density conditioning can create a trade-off. Scaling enlarges the range of $e^u$ and supplies no uniform bound on density conditioning.

\section{Attainment and the cost of approaching the boundary}\label{sec:endpoint}
Assume~\eqref{eq:base}. Local derivative norms use a fixed chart at $p$, and the reference form $m$ is held fixed. The covector $db(p)$ is the differential of $b$ at $p$, its vanishing means that all first coordinate derivatives vanish.
\begin{theorem}\label{th:endpoint}
A smooth compatible potential attains $\lambda_*=-b(p)$ with a positive Lyapunov rate if and only if
\begin{equation}\label{eq:criterion}
 db(p)=0.
\end{equation}
If this holds, every $V$ in~\eqref{eq:base} attains the endpoint after a sufficiently large scaling, retaining $c_0$.
There is $L_b<\infty$ such that $|b(x)-b(p)|\le L_bd$ on $K$. For $0<\eta<\lambda_*$, $sV$ is compatible with $(\lambda_*-\eta,c_0)$ whenever
\begin{equation}\label{eq:C2upper}
 s>0,\qquad s\ge\frac{L_b^2}{4c_0a_0\eta}.
\end{equation}
If $db(p)\ne0$, every compatible potential with volume rate $\lambda_*-\eta$ satisfies, on a fixed small coordinate ball $B$,
\begin{equation}\label{eq:C2lower}
 \|D^2u\|_{L^\infty(B)}\ge C/\eta
\end{equation}
for all small $\eta>0$. Thus the infimal $C^2$ cost has order $\eta^{-1}$.
Nevertheless, a smooth $q$ supported arbitrarily close to $p$ can make the divergence relative to $e^qm$ satisfy $d(b+f(q))(p)=0$, without changing $\lambda_*$.
\end{theorem}
\begin{proof}
Compatibility implies $du(p)=0$, hence $d(f(u))(p)=0$ because $f(p)=0$. At the endpoint, $b+f(u)\le b(p)$ with equality at the interior point $p$, so $db(p)=0$. Conversely, this condition gives $|b-b(p)|\le Cd^2$ near $p$. Since $-f(V)\ge c_0a_0d^2$, the ratio $(b-b(p))_+/[-f(V)]$ is bounded near $p$ and on the compact remainder. Large $s$ gives $b+sf(V)\le b(p)$ and preserves $c_0$.

Smoothness and compactness give $L_b$. For $r=d(x,p)>0$,
\[
 \frac{(b(x)-b(p)-\eta)_+}{-f(V)(x)}
 \le\frac{(L_br-\eta)_+}{c_0a_0r^2}
 \le\frac{L_b^2}{4c_0a_0\eta}.
\]
This proves~\eqref{eq:C2upper} and an $O(\eta^{-1})$ bound on the full $C^2$ norm for fixed $V$.
For the lower bound take $p=0$ and a unit vector $e$ with $db(0)e=\kappa>0$. On a small ball, $|f(x)|\le L|x|$ and $b(re)-b(0)\ge\kappa r/2$. At $r=4\eta/\kappa$, the volume inequality requires $f(u)(re)\le-\eta$. If $H=\|D^2u\|_{L^\infty(B)}$, then $|du(re)|\le Hr$, giving $\eta\le LHr^2$ and $H\ge\kappa^2/(16L\eta)$.
Finally $d(f(q))(p)=dq(p)A_f$, with $A_f$ invertible. Prescribe $dq(p)=-db(p)A_f^{-1}$ using a linear function times a local cutoff. Proposition~\ref{th:dual} preserves the value. This changes the reference form to $e^qm$, $q$ itself is not asserted to be a Lyapunov certificate. Relative to the old form the resulting logarithm is $q+sV$ and need not have a minimum at $p$, so this construction does not establish endpoint compatibility relative to the original $m$.
\end{proof}
The value is intrinsic to the flow on $K$, whereas Lyapunov compatibility of a density's logarithm is relative to the reference form. Growth of the potential Hessian reflects this local obstruction, independently of a chosen numerical algorithm. It is a statement about derivatives of the logarithmic density, not the Riemann curvature tensor, the latter vanishes identically in dimension one.

\emph{Pointwise attainment and a volume prefactor.}
Under~\eqref{eq:base}, the original reference volume attains the limiting exponent in an integral estimate, whether or not $db(p)=0$. Indeed, $V(\varphi_sx)\le e^{-c_0s}V(x)$ and the quadratic lower bound give
\[
 |b(\varphi_sx)-b(p)|\le\frac{L_b}{\sqrt{a_0}}\sqrt{V(x)}e^{-c_0s/2}.
\]
Writing $V_{\max}=\max_KV$, Liouville's formula~\eqref{eq:jacobian} therefore yields
\[
 J_m\varphi_t(x)\le
 \exp\!\left(\frac{2L_b\sqrt{V(x)}}{c_0\sqrt{a_0}}
 (1-e^{-c_0t/2})\right)e^{-\lambda_*t}.
\]
Consequently, for every measurable $E\subset K$,
\begin{equation}\label{eq:reference_endpoint}
 m(\varphi_tE)\le M_{\rm vol}e^{-\lambda_*t}m(E),\quad
 M_{\rm vol}=\exp\!\left(\frac{2L_b\sqrt{V_{\max}}}{c_0\sqrt{a_0}}\right),\quad t\ge0,
\end{equation}
where $m(E)=\int_Em$ and $J_m\varphi_t$ is the positive Jacobian relative to $m$. Estimate~\eqref{eq:reference_endpoint} establishes the limiting volume exponent with an explicit uniform prefactor. Together with the attainment criterion and derivative bounds, it separates integral contraction at the limiting exponent from realization by a compatible instantaneous certificate $D_u\le-\lambda_*$, which remains subject to the proved regularity obstruction.

\section{Sharp Sobolev costs and concentration at the limiting rate}\label{sec:sobolev}
The supremum norm detects the height of a derivative concentration but not its spatial extent. We now determine the sharp supercritical Sobolev cost in arbitrary dimension. The letter $r$ denotes the integrability exponent. The space $W^{2,r}(\Omega)$ consists of functions whose weak coordinate derivatives through order two lie in $L^r(\Omega)$, with the usual sum norm, for $r=\infty$ use essential suprema. Norms are unweighted coordinate norms on the stated domain, smooth background metrics give equivalent norms on its compact closure. For $n<r<\infty$, the Morrey embedding $W^{2,r}\subset C^{1,1-n/r}$ makes $du(p)$ meaningful~\cite{AdamsFournier2003}, $W^{2,\infty}$ has a $C^{1,1}$ representative.

\begin{definition}\label{def:sobolev}
Let $\Omega$ be a bounded smooth coordinate domain, $p\in\Omega$, with $K=\overline\Omega$ positively invariant. For $r>n$, including $r=\infty$, a function $u\in W^{2,r}(\Omega)$ is Sobolev compatible with $(\lambda,c)$ if its $C^1$ representative satisfies $u(p)=0$ and $ad^2\le u\le Ad^2$ on $\overline\Omega$ for some $a,A>0$, and $du[f]\le-cu$, $b+du[f]\le-\lambda$ almost everywhere in $\Omega$. The same definition is used for $n=r=1$, since $W^{2,1}$ functions on an interval have a $C^1$ representative. In these cases the differential inequalities extend to the closure by continuity. No representative-independent point condition is imposed for other exponents $r\le n$.
\end{definition}
For smooth costs we retain the original compact $K$, without requiring its boundary to be smooth. Let $\Omega$ be any bounded smooth coordinate domain with $p\in\Omega$ and $\overline\Omega\subset K$, equality $\overline\Omega=K$ is allowed. Define
\[
 \mathcal A_\eta(c)=\{u\in C^\infty(U):u\text{ satisfies~\eqref{eq:compatible}
 on }K\text{ with }(\lambda_*-\eta,c)\},
 \qquad 0<\eta<\lambda_*,
\]
where $a,A$ may vary with $u$, and set
\begin{equation}\label{eq:sobolevcost}
 \mathfrak C_r(\eta;c,\Omega)=
 \inf_{u\in\mathcal A_\eta(c)}\|u\|_{W^{2,r}(\Omega)}.
\end{equation}
An infimum over an empty class is $+\infty$. The domain $\Omega$ measures the cost near the equilibrium, while compatibility is still imposed on all of $K$.

In the following theorem all domains, charts and dynamical data, including the baseline $V$ and its comparison constants, are fixed. Constants in the uniform estimates may depend on these data and on $r$, and in fixed-$c$ cost estimates also on that prescribed $c$. They are independent of $\eta$, $R$ and the candidate potential $u$, including its freely varying quadratic constants $a,A$. The local obstruction in part~\textup{(1)} is also independent of the candidate's Lyapunov decay rate.

\begin{theorem}\label{th:sobolev}
Assume~\eqref{eq:base} and $db(p)\ne0$.
\begin{enumerate}
\item \textbf{Local obstruction.} For each $n<r<\infty$ and a fixed small coordinate ball $B$ at $p$, every smooth compatible potential with volume rate $\lambda_*-\eta$ satisfies
\begin{equation}\label{eq:morreycost}
 \|D^2u\|_{L^r(B)}\ge k_r\eta^{n/r-1}
\end{equation}
for all sufficiently small $\eta>0$, with $k_r>0$ independent of $u,\eta$. For $r=\infty$ this is~\eqref{eq:C2lower}. The same local lower bounds hold for $u\in W^{2,r}(B)$ satisfying quadratic bounds and the volume inequality almost everywhere. No nonnegative $C^1$ function with $u(p)=0$ satisfies the endpoint volume inequality near $p$. In particular there is no Sobolev-compatible endpoint in the classes of Definition~\ref{def:sobolev}.
\item \textbf{Sharp multidimensional cost.} Put
\begin{equation}\label{eq:rootparameters}
 B_0=\frac{L_b}{\sqrt{a_0}},\qquad C=\frac{4B_0}{c_0},
 \qquad\varepsilon=\frac{\eta}{B_0}.
\end{equation}
The smooth potentials, defined near $K$ and extendible to $U$,
\begin{equation}\label{eq:rootsmoothing}
 U_\eta=C\bigl(\sqrt{V+\varepsilon^2}-\varepsilon\bigr)
\end{equation}
are compatible with $(\lambda_*-\eta,c_0/2)$. For every fixed $0<c\le c_0/2$ and $n<r\le\infty$, there are $0<k_r\le K_r<\infty$ such that
\begin{equation}\label{eq:sharpSobolevND}
 k_r\eta^{n/r-1}\le\mathfrak C_r(\eta;c,\Omega)
 \le K_r\eta^{n/r-1}
\end{equation}
for small $\eta$, with $n/\infty=0$. Moreover $\|U_\eta\|_{C^1(K)}$ is uniformly bounded. The densities $e^{U_\eta}$ are uniformly bounded above, bounded below by one, and have uniformly bounded first derivatives on $K$.
For the larger range $0<c\le c_0$, scaling still gives the valid upper bound $\mathfrak C_r(\eta;c,\Omega)=O(\eta^{-1})$, its sharpness for $c>c_0/2$ is not asserted.
\item \textbf{Norm-budget deficit.} There is $\eta_0>0$, depending only on the fixed data and $r$, such that a budget $\|u\|_{W^{2,r}(B)}\le R$, $n<r<\infty$, forces, whenever $0<\lambda_*-\lambda<\eta_0$,
\begin{equation}\label{eq:budget}
 \lambda_*-\lambda\ge k'_r R^{-r/(r-n)}.
\end{equation}
For $r=\infty$ the lower bound is $k'/R$. These exponents are sharp in the freely conditioned smooth class. For fixed $0<c\le c_0/2$, let $\Lambda_r(R,c;\Omega)$ be the supremum of positive $\lambda$ for which some smooth potential is compatible on $K$ with $(\lambda,c)$ and satisfies $\|u\|_{W^{2,r}(\Omega)}\le R$, take zero for an empty class. For all sufficiently large $R$,
\begin{equation}\label{eq:budgetsharp}
 \lambda_*-\Lambda_r(R,c;\Omega)\asymp R^{-1/(1-n/r)},
 \qquad n<r\le\infty.
\end{equation}
The notation $\asymp$ means two inequalities with positive constants independent of the limiting parameter.
\item \textbf{Scalar borderline concentration.} In the scalar model
\begin{equation}\label{eq:scalarmodel}
 K=[-\ell,\ell],\quad f(x)=-x+ax^2,\quad m=\dd x,
 \qquad a>0,\quad0<a\ell<1,
\end{equation}
$\lambda_*=1$. Fix $0<c\le1-a\ell$ and abbreviate the cost on $\Omega=(-\ell,\ell)$ to $\mathfrak C_r(\eta)$. The following cost estimates hold for all sufficiently small $\eta>0$. For $1<r\le\infty$,
\begin{equation}\label{eq:sharpSobolev}
 k_r\eta^{1/r-1}\le\mathfrak C_r(\eta)\le K_r\eta^{1/r-1},
\end{equation}
whereas
\begin{equation}\label{eq:critical}
 0<k_1\le\mathfrak C_1(\eta)\le K_1<\infty.
\end{equation}
The explicit family below realizing the upper bounds converges uniformly and strongly in $W^{1,r}$ for every $1\le r<\infty$ to the Lipschitz function $u_0=C|x|$, with $C$ specified in the proof. This limiting object lies outside the original admissible class because it loses the quadratic upper bound at zero. Its second-derivative measures converge to a nonzero atom at zero. No $W^{2,1}$ Sobolev-compatible potential attains the endpoint. The stated atom describes this constructed family, not a classification of all optimizing sequences.
\end{enumerate}
\end{theorem}
\begin{proof}
\emph{Local lower bound and nonattainment.} Use the chart and $e,L,\kappa$ from Theorem~\ref{th:endpoint}. For $h=4\eta/\kappa$, the volume inequality at $he$ gives $|du(he)|\ge\eta/(Lh)$. On $B_{2h}$ the scaled Morrey estimate yields
\[
 |du(he)-du(0)|\le C_r h^{1-n/r}\|D^2u\|_{L^r(B_{2h})}.
\]
Since $du(0)=0$ by quadratic bounds, this implies
$\|D^2u\|_{L^r(B)}\ge\eta/(C_r Lh^{2-n/r})=k_r\eta^{n/r-1}$.
The gradient of a $W^{2,r}$ function, $r>n$, is continuous, so its almost-everywhere volume inequality extends pointwise and the same argument applies. The $r=\infty$ bound was proved in Theorem~\ref{th:endpoint}.
For any nonnegative $C^1$ endpoint candidate with $u(0)=0$, $du(0)=0$ and $f(x)=O(|x|)$ give $f(u)(x)=o(|x|)$. Along $he$ with $db(0)e=\kappa>0$, the endpoint inequality would require $\kappa h+o(h)\le0$, a contradiction.

\emph{Multidimensional construction.} Since $db(p)\ne0$, $L_b$ and $B_0$ can be chosen positive. On $K$, write $z=\sqrt V$. Formula~\eqref{eq:rootsmoothing} has quadratic bounds because
\begin{equation}\label{eq:rootquad}
 U_\eta=\frac{CV}{\sqrt{V+\varepsilon^2}+\varepsilon},\qquad
 \frac{Ca_0d^2}{\sqrt{V_{\max}+\varepsilon^2}+\varepsilon}
 \le U_\eta\le\frac{CA_0}{2\varepsilon}d^2,
\end{equation}
where $V_{\max}=\max_KV$. Differentiation along the flow gives
\begin{equation}\label{eq:rootdecay}
 f(U_\eta)=\frac{Cf(V)}{2\sqrt{V+\varepsilon^2}}
 \le-\frac{c_0}{2}U_\eta,
\end{equation}
since $V/\sqrt{V+\varepsilon^2}\ge\sqrt{V+\varepsilon^2}-\varepsilon$.
The Lipschitz bound for $b$ and $d\le z/\sqrt{a_0}$ imply
\begin{equation}\label{eq:rootvolume}
 b-b(p)+f(U_\eta)
 \le B_0z-\frac{2B_0z^2}{\sqrt{z^2+\varepsilon^2}}\le\eta.
\end{equation}
Indeed, for $z\le\varepsilon$ the right side is at most $B_0\varepsilon=\eta$, while for $z\ge\varepsilon$ it is at most $(1-\sqrt2)B_0z\le0$.

Smoothness on a neighborhood of $K$ follows because $V$ has a strict quadratic minimum at the interior point $p$ and is positive on the compact remainder of $K$. Thus some neighborhood of $K$ has $V\ge0$. Multiplying~\eqref{eq:rootsmoothing} by a cutoff equal to one near $K$ extends it smoothly to $U$ for each $\eta$, without changing any inequality or norm used here.

We next bound derivatives on $\Omega$. In coordinates centered at $p$, $V\asymp |x|^2$, $DV=O(|x|)$ and $D^2V=O(1)$ near zero. The exact identities are
\begin{equation}\label{eq:rootHessian}
 DU_\eta=\frac{CDV}{2\sqrt{V+\varepsilon^2}},\qquad
 D^2U_\eta=\frac{CD^2V}{2\sqrt{V+\varepsilon^2}}
 -\frac{C\,DV\otimes DV}{4(V+\varepsilon^2)^{3/2}}.
\end{equation}
Consequently $|DU_\eta|\le C_1$ and
$|D^2U_\eta|\le C_2(|x|^2+\varepsilon^2)^{-1/2}$ on a fixed small ball. On the compact remainder of $\overline\Omega$ the derivatives are uniformly bounded because $V$ has a positive minimum. For finite $r>n$, radial integration gives
\begin{equation}\label{eq:radialcost}
 \int_{|x|<\rho}|D^2U_\eta|^r\dd x
 \le C_3\int_0^\rho\frac{t^{n-1}}{(t^2+\varepsilon^2)^{r/2}}\dd t
 \le C_4\varepsilon^{n-r}.
\end{equation}
The function and gradient are bounded independently of $\eta$, and the Hessian supremum is $O(\varepsilon^{-1})$. Since $\varepsilon=\eta/B_0$, these are the upper bounds in~\eqref{eq:sharpSobolevND}. A ball $B\Subset\Omega$ supplies the matching lower bounds. The same derivative argument on finitely many charts covering $K$ proves the uniform $C^1(K)$ bound and then the stated bounds on $e^{U_\eta}$. For $c\le c_0$ the alternative $sV$ from~\eqref{eq:C2upper} gives $O(\eta^{-1})$ in every stated Sobolev norm and retains the full rate $c_0$.

\emph{Norm budgets.} Rearranging~\eqref{eq:morreycost} proves~\eqref{eq:budget}. Let $\gamma_r=1-n/r>0$. The two cost bounds give $\mathfrak C_r(\eta;c,\Omega)\asymp\eta^{-\gamma_r}$ at $c\le c_0/2$. In particular, choosing $\eta=(K_r/R)^{1/\gamma_r}$ in the explicit upper construction gives a feasible rate for large $R$. Conversely, every feasible rate with small gap obeys $\eta\ge(k_r/R)^{1/\gamma_r}$. Gaps outside this small regime are bounded below by a fixed positive number and obey the required lower bound for large $R$ as well. Taking the supremum of feasible rates proves~\eqref{eq:budgetsharp}, without assuming that this smooth supremum is attained. For $r=\infty$, use $\gamma_r=1$.

\emph{Scalar bounds and concentration.} In~\eqref{eq:scalarmodel}, $b=-1+2ax$ and $V=x^2$ has $f(V)\le-2(1-a\ell)V$, so $\lambda_*=1$. For $0<\eta<a\ell$, at $x_\eta=\eta/a$ the volume inequality requires
$f(x_\eta)u'(x_\eta)\le-\eta$, hence $u'(x_\eta)\ge a/(1-\eta)$. Since $u'(0)=0$, the fundamental theorem of calculus and H\"older's inequality give
\begin{equation}\label{eq:scalarHolder}
 \frac{a}{1-\eta}\le\int_0^{x_\eta}|u''|\dd x
 \le x_\eta^{1-1/r}\|u''\|_{L^r(0,x_\eta)},
 \qquad1\le r\le\infty.
\end{equation}
These are all the scalar lower bounds, including the positive bound at $r=1$.

Set $\delta=1-a\ell$, $C=4a/\delta$ and $\varepsilon=\eta/(2a)$. The specialization of~\eqref{eq:rootsmoothing} is
\begin{equation}\label{eq:smoothing}
 \begin{split}
 u_\eta(x)&=C(\sqrt{x^2+\varepsilon^2}-\varepsilon),\\
 u_\eta'(x)&=\frac{Cx}{\sqrt{x^2+\varepsilon^2}},\qquad
 u_\eta''(x)=\frac{C\varepsilon^2}{(x^2+\varepsilon^2)^{3/2}}.
 \end{split}
\end{equation}
It has quadratic bounds by~\eqref{eq:rootquad}, and $xu_\eta'\ge u_\eta$ gives
\begin{equation}\label{eq:scalarDecay}
 f(u_\eta)=-(1-ax)xu_\eta'\le-\delta u_\eta.
\end{equation}
For $x\le0$, $b+f(u_\eta)+1\le0$, for $0\le x\le\varepsilon$ it is at most $2ax\le\eta$. For $x\ge\varepsilon$ it is at most
$2ax-C\delta x^2/\sqrt{x^2+\varepsilon^2}\le(2-2\sqrt2)ax\le0$.
Thus it is compatible with every stated $c$. Its function and first derivative are uniformly bounded, and for finite $r\ge1$,
\begin{equation}\label{eq:integralnorm}
 \|u_\eta''\|_{L^r}^r=C^r\varepsilon^{1-r}
 \int_{-\ell/\varepsilon}^{\ell/\varepsilon}(1+y^2)^{-3r/2}\dd y.
\end{equation}
The integral has a positive finite limit, and $\|u_\eta''\|_\infty=C/\varepsilon$. This proves~\eqref{eq:sharpSobolev}--\eqref{eq:critical}.

The limit is $u_0=C|x|$. The bound $|u_\eta-C|x||\le C\varepsilon$ gives uniform convergence, and dominated convergence of first derivatives gives strong $W^{1,r}$ convergence for $1\le r<\infty$. Away from zero,
$f(u_0)=-(1-ax)C|x|\le-\delta u_0$ and $b+f(u_0)\le-1$, since $C\delta=4a\ge2a$.
The positive measures $u_\eta''\dd x$ converge weak-* to $2C\delta_0$: distributional differentiation gives the limit, and their total masses tend to $2C$ with no mass escaping at the endpoints. A $W^{2,1}$ function on an interval is $C^1$, so the already proved obstruction excludes a Sobolev-compatible endpoint. The Lipschitz limit has first derivative of bounded variation and second derivative a measure, not an $L^1$ function.
\end{proof}
The theorem determines the supercritical cost in every dimension, with the speed range in~\eqref{eq:sharpSobolevND} kept explicit. For $c_0/2<c\le c_0$, the lower bound remains valid and scaling the prescribed $V$ gives an $O(\eta^{-1})$ upper bound, a matching order is not asserted in this range. For $c>c_0$, this prescribed baseline provides no feasibility guarantee, although another Lyapunov function may allow a larger rate. The scalar bounded $W^{2,1}$ family is not uniformly integrable in its second derivatives. Its limit loses the quadratic upper bound and belongs to a different admissible class. Thus neither the bounded borderline mass nor its atomic limit contradicts the smooth or Sobolev nonattainment statements.

\section{Constrained searches and compactness}\label{sec:computation}
For the unrestricted divergence search, the final sign test in Appendix~\ref{app:perturbation} supplies a sufficient static criterion even without~\eqref{eq:base}. It does not by itself guarantee Lyapunov compatibility. We now impose that additional constraint in finite-dimensional and Sobolev classes.
Choose smooth basis functions $\psi=(\psi_1,\ldots,\psi_N)^\top$, write $u_\theta=\theta^\top\psi$ and $L_f\psi=(f(\psi_i))_{i=1}^N$, and let $q_0\in C(K)$ satisfy $m_0d^2\le q_0\le M_0d^2$ with $m_0,M_0>0$. For fixed $c,a,A>0$, consider
\begin{equation}\label{eq:SIP}
 \begin{aligned}
 \sup_{\theta,\lambda}\quad&\lambda\\
 \text{subject to}\quad&\lambda\ge0,\quad u_\theta(p)=0,\\
 &b+\theta^\top L_f\psi\le-\lambda,\\
 &\theta^\top L_f\psi+c\theta^\top\psi\le0,\\
 &a q_0\le\theta^\top\psi\le A q_0\qquad(x\in K).
 \end{aligned}
\end{equation}
\begin{proposition}\label{prop:computation}
Problem~\eqref{eq:SIP} is a semi-infinite linear program. If feasible, it attains a maximum no greater than $\lambda_*$. The parameter vector need not be unique, every positive feasible margin gives~\eqref{eq:distance} and~\eqref{eq:volume}.
For a Sobolev version let $K=\overline\Omega\Subset U$, where $\Omega$ is a bounded smooth coordinate domain containing $p$ and $K$ is positively invariant. Fix $n<r<\infty$, $r>1$, and $a,A,c,R>0$. Maximize $\lambda\ge0$ over $u\in W^{2,r}(\Omega)$ subject to
\begin{equation}\label{eq:weakprogram}
 \begin{split}
 &u(p)=0,\quad ad^2\le u\le Ad^2,\quad f(u)\le-cu,\\
 &b+f(u)\le-\lambda,\qquad\|u\|_{W^{2,r}(\Omega)}\le R.
 \end{split}
\end{equation}
If feasible, this convex problem attains a maximum. If $db(p)\ne0$, its maximum is strictly below $\lambda_*$ and obeys~\eqref{eq:budget} in the near-boundary regime.
\end{proposition}
\begin{proof}
All constraints in~\eqref{eq:SIP} are affine. Let $F=\{u_\theta|_K\}\subset C(K)$. If $u|_K=0$, positive invariance implies $u(\varphi_tx)=0$ and hence $f(u)|_K=0$. Thus $T:F\to C(K)$, $T(u|_K)=f(u)|_K$, is well defined and continuous by finite dimensionality. The order bounds make feasible functions bounded in $F$, and evaluation at $p$ bounds $0\le\lambda\le-b(p)$. The feasible set in $F\times\R$ is closed and bounded, hence compact.
For~\eqref{eq:weakprogram}, take a weakly convergent subsequence of a maximizing sequence in reflexive $W^{2,r}$. A further subsequence converges in $C^1(\overline\Omega)$ by compact embedding into $C^{1,\gamma}$ for $0<\gamma<1-n/r$. The inequalities and point condition pass to the limit. The norm constraint follows by weak lower semicontinuity. Linearity of $f(u)$ and convexity of the norm give convexity. Theorem~\ref{th:sobolev} proves the strict deficit.
\end{proof}
Joint optimization of $c$ and $\theta$ introduces a bilinear term. Fixed-$c$ slices remain linear. Feasibility is monotone in $c$, permitting an outer scalar search. Nested bases cannot shrink the feasible region, but convergence to the unrestricted rectangle needs approximation assumptions and removal of obstructive bounds. Polynomial data on semialgebraic $K$ admit sufficient sum-of-squares (SOS) relaxations implemented by semidefinite programming (SDP)~\cite{PrajnaParriloRantzer2004}. The regularity lower bounds give algorithm-independent limits on the achievable contraction margin under a prescribed norm budget, including for discretized searches and enlarged bases.
The matching norm-budget law~\eqref{eq:budgetsharp} holds in the freely conditioned smooth class at $0<c\le c_0/2$. Its scalar specialization is $\eta\asymp R^{-r/(r-1)}$. Fixed constants $a,A$ in~\eqref{eq:weakprogram} can impose additional obstructions, so the matching upper law is not asserted for that more restricted problem or for a particular sequence of basis searches.

\section{Examples and interpretation}\label{sec:examples}
\begin{figure}[h]
\center{\includegraphics[width=0.7\linewidth]{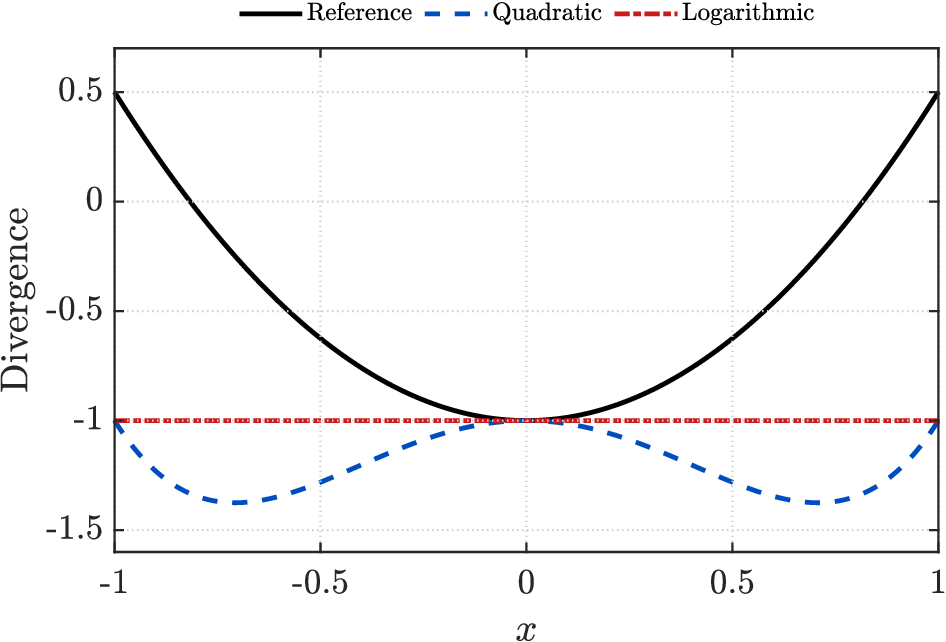}}
\caption{Both compatible potentials in Example~\ref{ex:cubic} attain the optimal uniform margin. The logarithmic potential makes divergence constant.}
\label{fig:cubic}
\end{figure}

\begin{example}\label{ex:cubic}
On $K=[-1,1]$ take $f=-x+\tfrac12x^3$ and $b=-1+\tfrac32x^2$. The field points inward and $V=x^2$ has $f(V)=-2x^2+x^4\le-V$. Hence $\I_f(K)=\{\delta_0\}$ and $\lambda_*=1$, although $b>0$ for $|x|>\sqrt{2/3}$. Since $b'(0)=0$, the endpoint is compatible. The potential $u=\tfrac32x^2$ gives
\[
 f(u)=-3x^2+\tfrac32x^4\le-u,\qquad
 D_u=-1-\tfrac32x^2+\tfrac32x^4\le-1.
\]
Alternatively,
\begin{equation}\label{eq:logpotential}
 u_*=-\tfrac32\log(1-x^2/2),\qquad
 e^{u_*}=(1-x^2/2)^{-3/2}
\end{equation}
is smooth on a neighborhood of $K=[-1,1]$, since $1-x^2/2\ge1/2$ on $K$, and gives $D_{u_*}\equiv-1$. It is compatible: $f(u_*)=-3x^2/2$ and $3x^2/4\le u_*\le3x^2/2$, using the elementary bounds on $-\log(1-z)$ for $0\le z\le1/2$. Thus $f(u_*)\le-u_*$. Figure~\ref{fig:cubic} compares the exact divergences.
\end{example}

\begin{figure}[h]
\center{\includegraphics[width=0.7\linewidth]{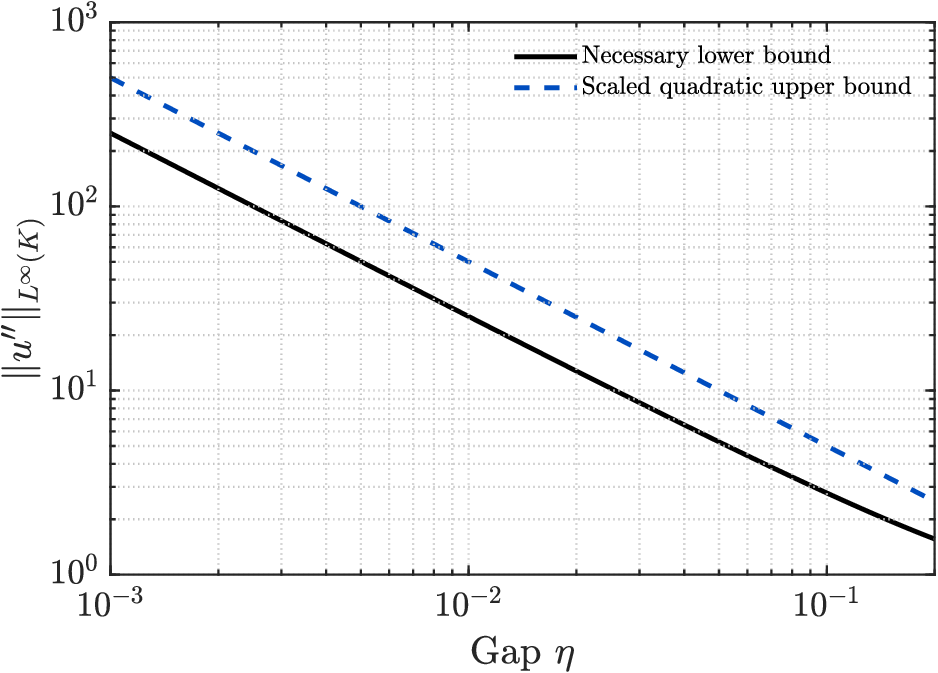}}
\caption{Bounds in Example~\ref{ex:scalar} for $a=1/2$, $\ell=1$. Both grow as $\eta^{-1}$. Theorem~\ref{th:sobolev} resolves the integral-norm costs}\label{fig:cost}
\end{figure}

\begin{example}\label{ex:scalar}
For~\eqref{eq:scalarmodel}, the unrestricted smooth $q=-2\log(1-ax)$ gives $b+f(q)=-1$. But $q'(0)=2a\ne0$, so it is not a positive-definite Lyapunov function. Besides~\eqref{eq:smoothing}, scaled quadratics
\[
 u_\eta=s_\eta x^2,\qquad s_\eta=\frac{a^2}{2(1-a\ell)\eta}
\]
give $D_{u_\eta}+1\le2ax-2s_\eta(1-a\ell)x^2\le\eta$ and preserve the Lyapunov rate $2(1-a\ell)$. Every compatible potential satisfies, by~\eqref{eq:scalarHolder},
\begin{equation}\label{eq:scalarC2}
 \|u''\|_\infty\ge\frac{a^2}{\eta(1-\eta)},\qquad0<\eta<a\ell,
\end{equation}
while the quadratic upper bound is $a^2/[(1-a\ell)\eta]$. Both have the sharp $C^2$ order (Fig.~\ref{fig:cost}). The localized family~\eqref{eq:smoothing} also resolves the concentration and gives sharp finite-$r$ Sobolev bounds. The same supremum-norm order need not imply the same integral-norm behavior.
\end{example}

\begin{example}\label{ex:multidimensional}
Let $n\ge2$, $0<a<1$, and $K=\{x\in\R^n:|x|\le1\}$. Consider
\[
 f_1=-x_1+ax_1^2,\qquad f_j=-2x_j\quad(2\le j\le n),
 \qquad m=\dd x_1\cdots\dd x_n.
\]
For $V=|x|^2$, $f(V)\le-2(1-a)V$, so the field points strictly inward on the boundary and~\eqref{eq:base} holds with $a_0=A_0=1$, $c_0=2(1-a)$. Here $b=-(2n-1)+2ax_1$, $\lambda_*=2n-1$ and $db(0)\ne0$. Theorem~\ref{th:sobolev} gives
\[
 U_\eta(x)=\frac{4a}{1-a}\bigl(\sqrt{|x|^2+\varepsilon^2}-\varepsilon\bigr),
 \qquad\varepsilon=\frac{\eta}{2a},
\]
with volume rate $2n-1-\eta$ and Lyapunov rate $1-a$. With $C=4a/(1-a)$, its Hessian has radial eigenvalue $C\varepsilon^2/(|x|^2+\varepsilon^2)^{3/2}$ and tangential eigenvalue $C/(|x|^2+\varepsilon^2)^{1/2}$ of multiplicity $n-1$. Therefore for every $r>n$ its $W^{2,r}$ norm has precisely the optimal order $\eta^{n/r-1}$. The density and its first derivatives remain uniformly bounded as $\eta\downarrow0$, while the Hessian cost diverges. For the critical exponent $r=n$, this family has Hessian norm of order $(\log(1/\eta))^{1/n}$, this last computation describes the family and is not claimed to be a sharp lower bound for all admissible potentials.
\end{example}

\begin{example}\label{ex:coupledoscillator}
Consider the coupled nonlinear oscillator
\[
 \dot x=y,\qquad \dot y=-x-y+\tfrac14xy,
\]
equivalently $\ddot x+(1-x/4)\dot x+x=0$, describing
a unit-stiffness oscillator with position-dependent damping.
Set
\[
 V=x^2+xy+y^2,\qquad K=\{V\le1\},\qquad m=dx\,dy.
\]
Since $\frac12|(x,y)|^2\le V\le\frac32|(x,y)|^2$,
$|xy|\le V$ and $|x+2y|\le2\sqrt V$, one obtains
\[
 L_fV=-V+\tfrac14xy(x+2y)
 \le-V+\tfrac12V^{3/2}\le-\tfrac12V
 ~~\text{on }K.
\]
Thus $K$ is positively invariant and~\eqref{eq:base} holds
with $a_0=1/2$, $A_0=3/2$, $c_0=1/2$.
Moreover,
\[
 b=-1+x/4,\qquad \lambda_*=1,\qquad
 db(0)=(1/4,0)\ne0.
\]
The reference divergence gives only the uniform margin
$1-1/(2\sqrt3)$, whereas adapted densities approach $1$.
Theorem~\ref{th:endpoint} excludes smooth compatible
attainment of this limiting margin.

Both approximation mechanisms are explicit. For $0<\eta<1$,
\[
 u_\eta^{\mathrm q}=\frac{V}{24\eta},\qquad
 L_fu_\eta^{\mathrm q}\le-\tfrac12u_\eta^{\mathrm q},
 \qquad b+L_fu_\eta^{\mathrm q}\le-1+\eta,
\]
where the last inequality follows from
$x/4\le\sqrt V/(2\sqrt3)$ and maximizing the resulting
quadratic in $\sqrt V$. Conversely, at $(8\eta,0)$,
any compatible potential with margin $1-\eta$ satisfies
$u_y(8\eta,0)\ge1/8$. Since $Du(0)=0$, this gives
\[
 \|D^2u\|_{L^\infty(K)}\ge\frac1{64\eta},
 \qquad 0<\eta<1/8,
\]
using the Euclidean operator norm.

The regularized family
$ U_\eta=\frac4{\sqrt3}
 \left(\sqrt{V+12\eta^2}-2\sqrt3\,\eta\right)$
instead satisfies
\[
 L_fU_\eta\le-\tfrac14U_\eta,\qquad
 b+L_fU_\eta\le-1+\eta,\qquad
 0\le U_\eta\le4/\sqrt3,\quad |DU_\eta|\le2\sqrt2.
\]
Indeed, these follow from the construction in
Theorem~\ref{th:sobolev}, using
$b-b(0)\le\sqrt V/(2\sqrt3)$.
For $\Omega=\operatorname{int}K$, fixed $0<c\le1/4$
and $2<r\le\infty$, that theorem yields
$\mathfrak C_r(\eta;c,\Omega)\asymp\eta^{2/r-1}$ and 
$ 1-\Lambda_r(R,c;\Omega)\asymp R^{-r/(r-2)}$,
with exponent $1$ in the budget law when $r=\infty$.
Thus the coupled oscillator exhibits both preservation
of the baseline Lyapunov rate by scaling and sharp
Sobolev costs with uniformly bounded density and gradient.
\end{example}

\begin{example}\label{ex:cycle}
On the Euclidean annulus $1\le r\le2$, put $a=3/2$ and $f=-r^2(r-a)\partial_r+r^{-1}\partial_\theta$. The divergence is $b=-4r^2+3ar$. Both boundaries point inward. Since $W=(r-a)^2$ satisfies $f(W)=-2r^2(r-a)^2$, invariant probabilities are supported on $r=a$. Motion on this cycle is uniform, giving a unique invariant probability and $\beta_K=-a^2=-9/4$. The radial potential $u_*=-4\log r+a/r$ gives
\[
 f(u_*)=(r-a)(4r+a)=4r^2-3ar-a^2,\qquad D_{u_*}=-a^2.
\]
The volume rate $9/4$ is attained, while trajectories approach a cycle. Its tangential exponent is zero and its transverse exponent is $-9/4$. The equilibrium compatibility hypotheses do not hold.
\end{example}
\begin{example}
For $\dot x=x$, $\dot y=-3y$, divergence is $-2$ but the origin is unstable. No quadratic decaying Lyapunov function exists on a full neighborhood. Negative divergence cannot replace the stability hypothesis.
\end{example}

\section{Concluding remarks}
The results characterize optimal weighted volume contraction on a prescribed positively invariant region. Under the quadratic Lyapunov hypothesis, the equilibrium trace determines the optimal value, while the nonlinear regional constraints govern its compatible realization, boundary attainment and regularity cost. Scaling preserves the full baseline decay rate. The fixed-reference first-order divergence defect gives the exact smooth pointwise attainment criterion and sharp inverse-gap Hessian cost. Estimate~\eqref{eq:reference_endpoint} additionally realizes the limiting volume exponent in the reference measure with an explicit uniform prefactor, including when the compatible instantaneous bound is unattainable.

For $r>n$, a regularized square root of the baseline Lyapunov function attains the sharp multidimensional $W^{2,r}$ cost at each fixed $0<c\le c_0/2$. It also gives the matching norm-budget deficit while keeping the density and its first derivatives uniformly bounded. The scalar borderline construction has bounded second-derivative mass and a Lipschitz limit with atomic second derivative. That limit leaves the original quadratically bounded class. Further questions concern sharp costs at larger prescribed Lyapunov rates and at critical integrability in higher dimensions. For the established supercritical range, the matching bounds and explicit construction give a complete cost law for approaching the optimal divergence bound.

\appendix
\section{Transport and conditional stability consequences}\label{app:transport}
The following auxiliary facts retain the geometric and stability consequences underlying the optimality results.
\subsection{Weights and evolving geometries}
For a background metric $g$, $dV_g$ is its volume form, $\nabla_g$ its gradient, $\Hess_g$ its covariant Hessian, and $\Delta_g=\diver_g\nabla_g$ its Laplace--Beltrami operator. For a covariant tensor $H$, $\tr_gH=g^{ij}H_{ij}$, summing repeated indices, $(g^{ij})$ is the inverse metric matrix. For positive $w$ and $\alpha\in\R$, define
\[
 D_{w^\alpha}f=\diver_gf+\alpha f(\log w)
 =\frac{\diver_g(w^\alpha f)}{w^\alpha}.
\]
The metric $\widetilde g_\alpha=w^{2\alpha/n}g$ has volume $w^\alpha dV_g$ and divergence $D_{w^\alpha}f$. For $w=e^{n\phi}$ it is $e^{2\alpha\phi}g$. For $\alpha=0$, both are independent of $w$. If $f=-\nabla_gV$, the divergence inequality reads
$-\Delta_gV-\alpha n\langle\nabla_gV,\nabla_g\phi\rangle_g\le-\lambda$.
For $\alpha\ne0$, $w^\alpha$ ranges over all positive smooth densities. Thus the classes admitting a contracting positive density, a conformal metric and an arbitrary metric coincide. Their complement contains the flows admitting none on the prescribed region.
If $w\ge0$ vanishes, logarithmic divergence is defined only where $w>0$. For smooth nonnegative $w$ and $\alpha\ge1$, the product
\[
 w^\alpha D_{w^\alpha}f=w^\alpha\diver_gf+\alpha w^{\alpha-1}f(w)
\]
extends continuously to the zero set (the second term is $f(w)$ for $\alpha=1$). This does not define a finite logarithmic divergence at a zero.

Let $g(t)$ be smooth, $w(x,t)>0$ be $C^1$, and $D(t)=\varphi_t(D)$ for a region with piecewise smooth boundary. On any finite interval of existence, differentiation of $\int_D\varphi_t^*(w^\alpha dV_{g(t)})$ gives
\begin{equation}\label{eq:transport}
 \begin{split}
 \frac{d}{dt}\int_{D(t)}w^\alpha dV_{g(t)}
 =\int_{D(t)}\bigl(D_{w^\alpha}f+\alpha\partial_t\log w +\tfrac12\tr_{g(t)}\partial_tg\bigr)w^\alpha dV_{g(t)}.
 \end{split}
\end{equation}
Indeed, $\partial_tdV_g=\tfrac12\tr_g(\partial_tg)dV_g$ and $\mathcal L_fdV_g=(\diver_gf)dV_g$, while the two derivatives of $w^\alpha$ give the other terms. If the parenthesized factor is at most $-\mu<0$ on a positively invariant region, integration gives volume decay with prefactor one. More generally, exponential weighted-volume dissipativity means
\[
 \int_{\varphi_tD}w(x,t)^\alpha dV_{g(t)}
 \le C e^{-\lambda t}\int_Dw(x,0)^\alpha dV_{g(0)}
\]
for all admissible regions and times, with $C,\lambda>0$. Allowing $\lambda=0$ gives the corresponding bounded-growth version. This use of dissipativity describes phase volume, rather than a storage and supply inequality~\cite{Willems72}. In the static case~\eqref{eq:transport} proves~\eqref{eq:volume}, first for regular regions and then by measure approximation. Strict negativity on compact $K$ extends to a neighborhood, which need not itself be positively invariant.
On a compact manifold without boundary, Stokes' theorem~\cite{Besse} gives
$\int_M D_{w^\alpha}f\,w^\alpha dV_g=0$.
Thus no smooth positive weight or metric gives uniform negative divergence on the entire manifold. This does not exclude proper invariant subsets.

\subsection{Stability from weights and gradient potentials}
Suppose $p$ is an isolated equilibrium, $w\in C^1$ is positive off $p$ and zero at $p$, and $V=w^\alpha$, $\alpha>0$, extends to $C^1$ across $p$. Off $p$,
\begin{equation}\label{eq:weightLyapunov}
 f(V)=V(D_{w^\alpha}f-\diver_gf).
\end{equation}
If $D_{w^\alpha}f\le-\gamma$, $\diver_gf\ge-d_0$, $\gamma>d_0\ge0$, and $c_1d^\nu\le V\le c_2d^\nu$ with $c_1,c_2,\nu>0$, then $f(V)\le-(\gamma-d_0)V$. Standard Lyapunov integration~\cite{Khalil2002} yields
\[
 d(\varphi_tx,p)\le(c_2/c_1)^{1/\nu}
 e^{-(\gamma-d_0)t/\nu}d(x,p)
\]
on a small positively invariant sublevel set. Such a set follows from the two-sided bounds and decay. A finite $d_0$ exists on compact neighborhoods by continuity. If evolving weights converge to $\widehat w$ satisfying these assumptions, apply the criterion to that fixed limit. Convergence alone supplies neither a relative inequality near its zero nor convergence of a geometric evolution to a valid certificate.

For a compact smooth invariant submanifold $N$ without boundary lying in the interior of the ambient manifold, squared distance $\rho_N^2$ is smooth on a small tubular neighborhood, where $\rho_N(x)=d_g(x,N)$. For positive $C^1$ $w$, bounded above and away from zero there, set $V_N=\rho_N^2w^\alpha$. Off $N$,
\[
 f(V_N)=V_N\left(D_{w^\alpha}f-\diver_gf+2\frac{f(\rho_N)}{\rho_N}\right).
\]
If $D_{w^\alpha}f+2f(\rho_N)/\rho_N\le-\gamma$ and $\diver_gf\ge-d_0$, with $\gamma>d_0$, this extends to $f(V_N)\le-(\gamma-d_0)V_N$ throughout the tube. Small sublevel tubes are invariant and
$\rho_N(\varphi_tx)\le C e^{-(\gamma-d_0)t/2}\rho_N(x)$.
The square is essential: an odd codimension power of distance need not be smooth across $N$.

If $f=-\nabla_{\hat g}V$ on a geodesically convex ball at $p$, $\nabla_{\hat g}V(p)=0$ and $\Hess_{\hat g}V\succeq\mu\hat g$, $\mu>0$, shrink the ball to get an upper Hessian bound $L\hat g$. Strong convexity gives
\[
 \tfrac\mu2 d_{\hat g}^2\le V-V(p)\le\tfrac L2 d_{\hat g}^2,
 \qquad\|\nabla_{\hat g}V\|^2\ge2\mu(V-V(p)).
\]
Thus $f(V)\le-2\mu(V-V(p))$ and
$d_{\hat g}(\varphi_tx,p)\le\sqrt{L/\mu}\,e^{-\mu t}d_{\hat g}(x,p)$
on a small invariant sublevel set. A limiting geometric flow can supply the metric. Stability still uses these additional gradient and convexity assumptions.

\section{Perturbing potentials and a static sign test}\label{app:perturbation}
The scaling proof makes an auxiliary-potential search unnecessary under~\eqref{eq:base}. The following construction is useful when preserving a selected nearly optimal potential.
\begin{lemma}\label{lem:jet}
For smooth $q$ and $\varepsilon>0$ there is $h\in C_c^\infty(U)$ supported arbitrarily close to $p$ with $h(p)=0$, $dh(p)=-dq(p)$ and $\|f(h)\|_{C(K)}<\varepsilon$.
\end{lemma}
\begin{proof}
In coordinates at $p$, choose a smooth cutoff $\chi$ equal to one on the unit ball and zero outside radius two. Put $h_\rho(x)=-\chi(x/\rho)dq(p)[x]$. Its gradient is bounded independently of small $\rho$, while $f=O(\rho)$ on its support. Hence $\|f(h_\rho)\|_\infty=O(\rho)$ and the prescribed differential holds.
\end{proof}
If $b+f(q)\le-\lambda_1$ and $\lambda<\lambda_1$, use an error below $\lambda_1-\lambda$ and set $\zeta=q-q(p)+h$. Then $\zeta(p)=0$, $d\zeta(p)=0$ and $b+f(\zeta)\le-\lambda$. Taylor estimates and compactness under~\eqref{eq:base} give $|\zeta|\le C_0V$, $|f(\zeta)|\le C_1V$. For large $s$,
\[
 \begin{split}
 u_s=\zeta+sV,\qquad &(s-C_0)V\le u_s\le(s+C_0)V,\\
 f(u_s)&\le-\frac{sc_0-C_1}{s+C_0}u_s.
 \end{split}
\]
The volume inequality is preserved by $f(V)\le0$, and the displayed decay rate tends to $c_0$. This retains the auxiliary-function proof of no loss and of the joint rectangle.

For fixed conditioning constants define $\lambda_{\rm LC}(K,c,a,A)$ as the supremum of positive compatible margins, taking zero for an empty class. Then
$\lambda_{\rm LC}(K,c,a,A)\le\max\{0,\lambda_*\}$, and decreasing $c$ enlarges the class. Write $\widehat\lambda_{\rm LC}(K,f)$ for the supremum over all positive $c,a,A$, with the same zero convention. In general it is at most $\max\{0,\lambda_*\}$, and it equals $\lambda_*$ under~\eqref{eq:base}. The zero convention matters if a positive-margin class is empty.

A sufficient sign test is available without~\eqref{eq:base}. If $\phi$ is smooth near $K$, $f(\phi)\ge0$ on $K$, and $b<0$ on $Z=\{f(\phi)=0\}\cap K$, then some $u=-C\phi$, $C\ge0$, has $b+f(u)\le-\lambda<0$. Indeed, $b$ is strictly negative near compact $Z$. On the compact remainder, $f(\phi)$ has a positive minimum, so large $C$ works without worsening the sign near $Z$. If $Z$ is empty use its positive minimum on all of $K$. If the negative-divergence neighborhood covers $K$, take $C=0$. This is a static test, not a convergence theorem for geometric flows.

\backmatter
\section*{Statements and Declarations}
\subsection*{Competing interests}

\textit{The author declares no competing interests.}
\subsection*{Data availability}
No external datasets were used. The figures illustrate explicit
functions given in the text.

\begingroup
\interlinepenalty=10000

\endgroup
\end{document}